\documentclass[12pt,a4paper]{article}

\usepackage[dutch,english]{babel}
\usepackage[UKenglish]{isodate}
\usepackage{amssymb, amsmath, calligra, mathrsfs, color, fancyhdr, enumitem, amsthm, hyperref, graphicx, picins, array, tikz-cd, subfiles}
\usepackage[toc,page]{appendix}
\usepackage[a4paper, top=3.5cm, bottom=3cm, left=2.5cm, right=3.5cm]{geometry}
\usepackage[all]{xy}
\newtheorem{theorem}{Theorem}
\newtheorem{lemma}[theorem]{Lemma}
\newtheorem{proposition}[theorem]{Proposition}
\newtheorem{corollary}[theorem]{Corollary}
\newtheorem{conjecture}[theorem]{Conjecture}
\theoremstyle{definition}
\newtheorem{definition}[theorem]{Definition}
\newtheorem{remark}[theorem]{Remark}
\newtheorem{example}[theorem]{Example}

\renewcommand{\O}{\mathcal{O}}

\newcommand{\commentaar}[1]{}
\newcommand{\C}{\ensuremath{\mathbb{C}}}

\newcommand{\Q}{\ensuremath{\mathbb{Q}}}
\newcommand{\Z}{\ensuremath{\mathbb{Z}}}

\newcommand{\Hom}{\textrm{Hom}}
\DeclareMathOperator{\ShHom}{\mathscr{H}\text{\kern -3pt {\calligra\large om}}\,}
\renewcommand{\O}{\mathcal{O}}

\newcommand{\tors}{\mathrm{tors}}

\newcounter{nootje}
\renewcommand\check[1]{[*\thenootje]\marginpar{\tiny\begin{minipage}{25mm}\begin{flushleft}\thenootje : #1\end{flushleft}\end{minipage}}\addtocounter{nootje}{1}}

\DeclareFontFamily{U}{wncy}{}
\DeclareFontShape{U}{wncy}{m}{n}{<->wncyr10}{}
\DeclareSymbolFont{mcy}{U}{wncy}{m}{n}
\DeclareMathSymbol{\Sh}{\mathord}{mcy}{"58} 

\newcommand{\onlyinsubfile}[1]{#1}

\newcommand{\Qfive}{\ensuremath{\Q\!\left(\sqrt[4]5\,\right)}}
\newcommand{\Qfivei}{\ensuremath{\Q\!\left(i, \sqrt[4]5\,\right)}}
\newcommand{\Jac}{\ensuremath{\mathop\mathrm{Jac}}}
\newcommand{\Res}{\ensuremath{\mathop\mathrm{Res}}}

\begin{document}

\title{The Birch and Swinnerton-Dyer conjecture\\ for an elliptic curve over $\Qfive$}
\author{Raymond van Bommel}
\cleanlookdateon
\maketitle

\renewcommand{\onlyinsubfile}[1]{}

{\bf Abstract.} In this paper we show the Birch and Swinnerton-Dyer conjecture for a certain elliptic curve over $\Qfive$ is equivalent to the same conjecture for a certain pair of hyperelliptic curves of genus 2 over $\Q$. We numerically verify the conjecture for these hyperelliptic curves. Moreover, we explain the methods used to find this example, which turned out to be a bit more subtle than expected.


\section{Introduction}

The Birch and Swinnerton-Dyer conjecture (\cite{BirchSwinnertonDyer}) has been generalised by Tate (\cite{TateBourbaki}) to abelian varieties of higher dimension and over general number fields.

\begin{conjecture}[\textrm{BSD, \cite[Conj.\ 2.10, p.\ 224]{Gross}}]
Let $A/K$ be an abelian variety of dimension $d$ and algebraic rank $r$ over a number field $K$ of discriminant $\Delta$. Let $L(s)$ be its $L$-function, $A^{\vee}$ its dual, $R$ its regulator, $\Sh$ its Tate-Shafarevich group and $\Omega$ the product of its real and complex periods. For each prime $\mathfrak{p}$ of $\O_K$, let $c_{\mathfrak{p}}$ be the Tamagawa number of $A$ at $\mathfrak{p}$. Then $\Sh$ is finite, $L(s)$ admits an analytic continuation to $\C$ having a zero of order $r$ at $s = 1$, and $$\lim_{s \rightarrow 1} {(s-1)^{-r} L(s)} = \frac{\Omega \cdot R \cdot |\Sh| \cdot \prod_{\mathfrak{p}} c_{\mathfrak{p}}}{|A(K)_{\tors}| \cdot |A^\vee(K)_{\tors}| \cdot |\Delta|^{d/2}}.$$
\end{conjecture}

In 1989, Kolyvagin (\cite{Kolyvagin, KolyvaginSurvey}) proved equality of the analytic and algebraic rank for modular elliptic curves over $\Q$ of analytic rank at most 1. After the proof of the modularity theorem (\cite{Modularity}), this part of the conjecture is now known for all elliptic curves over $\Q$ of analytic rank at most 1.

For elliptic curves with complex multiplication more is known. In 1991, Rubin (\cite{Rubin}) proved the correctness of the $p$-part of BSD for elliptic curves over an imaginary quadratic field $K$ with complex multiplication by $K$, analytic rank equal to 0, and $p$ coprime to $|\O_K^*|$.

Originally, the Birch and Swinnerton-Dyer conjecture has been conceived based on numerical calculations with elliptic curves. In \cite{NumericalVerification}, the author numerically verified the conjecture for hundreds of hyperelliptic curves of genus 2 and 3 over $\Q$, extending the work of Flynn, Lepr\'evost, Schaefer, Stein, Stoll and Wetherell (\cite{EmpiricalEvidence}), who numerically verified BSD for 32 modular hyperelliptic curves of genus 2 over $\Q$, using modularity.

This verification consists of two parts. First, we check that analytic rank (established numerically) and the algebraic rank are equal. Then we numerically compute all terms in the BSD formula except for $|\Sh|$ (to more than 20 digits precision), and by rearranging the formula we deduce a predicted value for $|\Sh|$. This will a priori be some real number, but if the BSD conjecture is true then it should in fact be the square of a positive integer, cf. earlier results of Poonen and Stoll ([PoSt99]). So if our conjectural value of $|\Sh|$ is indeed the square of a positive integer to high precision, then this provides strong numerical evidence for the conjecture.


After finishing this verification, a natural question that arose was if the numerical verification for genus 2 curves over $\Q$, could provide us with examples of elliptic curves $E$ over quadratic number fields for which BSD numerically seems to hold. The Weil restriction of $E$ to $\Q$ is an abelian variety of dimension 2 over $\Q$ and might have the chance of being the Jacobian of a genus 2 curve over $\Q$. As the Jacobi locus is dense in the moduli space, one might expect this to happen very often. This was not the case. While trying many examples, all seemed to fail.

However, this Weil restriction becomes a product of two elliptic curves, after base change. The product of two elliptic curves, taken with the associated product polarisation, does not lie in the Jacobi locus. The best we could hope for is the existence of another polarisation, which makes it isomorphic (as polarised abelian variety) to the Jacobian of a curve of genus 2. This is actually only possible in a few special cases.
By trying other polarisations in these special cases, we found an example of an elliptic curve over $\Q(\sqrt5)$, whose Weil restriction is isogenous to the Jacobian of a curve of genus 2 over $\Q$. However, the isogeny was only defined over $\Q(\sqrt[8]{5},i)$. We applied some reduction steps to reduce the size of this field and arrive at the following theorem
\begin{theorem}\label{thm:BSDoverQ5}
Let $E$ over $\Qfive$ be the elliptic curve given by
$$y^2 = x^3 + \sqrt[4]{5} \cdot x^2 - \left(5 + 3\sqrt{5}\right) \cdot x + \sqrt[4]{5} \left( 5 + \sqrt{5}\right).$$
Let $H$ and $H'$ over $\Q$ be the hyperelliptic curves given by $y^2 = x^5 - x^3 + \tfrac15 \cdot x,$ and $y^2 = x^5 - 5 \cdot x^3 + 5\cdot x,$ respectively. Then the generalised Birch and Swinnerton-Dyer conjecture holds for $E$ over $\Qfive$ if and only if it holds for the Jacobians $\Jac{H}$ and $\Jac{H'}$ over $\Q$.
\end{theorem}
Finally, because of this reduction of the size of the field, we were able to numerically verify the BSD conjecture for the mentioned hyperelliptic curves.

We could also phrase the problem we solved as a moduli problem. For fixed $N$, we consider the space $\mathcal{M}$ of quintuples $(E_1, E_2, A, \phi, \rho)$, where $E_1$ and $E_2$ are elliptic curves, $(A, \phi)$ is a principally polarised abelian surface, and $\rho : E_1 \times E_2 \rightarrow A$ is an isogeny of degree $N$. If $\iota : \mathcal{M} \rightarrow \mathcal{M}$ is the involution that swaps $E_1$ and $E_2$, then our problem is the finding of rational points of $\mathcal{M} / \iota$, for which $(A, \varphi)$ is not a product of elliptic curves.

This moduli problem (or variations thereof) has been studied extensively by others. This started with Hayashida and Nishi in \cite{HayashidaNishi}. More recently, there is work of Rodriguez-Villegas (\cite{RodriguezVillegas}), Lange (\cite{Lange}), and Kani (\cite{Kani14}, \cite{Kani16}). However, as far as we are aware, none of these results gives a way to control the size of the field of definition for the isogeny $\rho$, which is needed for our verification of the BSD conjecture.

The organisation of this article is as follows. In the first section, the final results will be shown, the equivalence of BSD for a certain elliptic curve over a quartic field and BSD for a certain pair of hyperelliptic curves of genus 2 over $\Q$. In the second section, the methods used to find this example will be demonstrated. First we study which elliptic curves could have to potential to become isogenous to the Jacobian of a genus 2 curve after Weil restriction. Then we explain how the required isogenies, which are very easy to find analytically, were algebraised. Finally, we describe some steps that had to be taken to reduce the size of the number field over which these maps are defined, which was actually necessary to be able to complete the verification.

The author wishes to thank his supervisors David Holmes and Fabian Pazuki, and Maarten Derickx for useful discussions that led to improvements of this article.

\onlyinsubfile{

}

\section{Verification for an elliptic curve over $\Qfive$}

Throughout this section, let $E$ be the elliptic curve over $\Qfive$ given by the Weierstra\ss\ equation $$y^2 = x^3 + \sqrt[4]{5} \cdot x^2 - \left(5 + 3\sqrt{5}\right) \cdot x + \sqrt[4]{5} \left( 5 + \sqrt{5}\right).$$ Even though it has $j$-invariant $282880\sqrt{5} + 632000$, it is not the base change of an elliptic curve over $\Q(\sqrt5\,)$, which can be verified using the isomorphism criteria from \cite[Sect.\ III.1, p.\ 42--51]{Silverman}. The curve $E$ geometrically has complex multiplication by $\Z[\sqrt{-5}]$.

Let $H$ be the hyperelliptic curve of genus 2 over $\Q$ given by the Weierstra\ss\  equation $y^2 = x^5 - x^3 + \tfrac15 \cdot x.$ Let $H' \colon y^2 = x^5 - 5 \cdot x^3 + 5\cdot x$ over $\Q$ be the quadratic twist of $H$ over $\Q(\sqrt{5}\,)$.

\onlyinsubfile{{\color{gray}
In order to prove Theorem \ref{thm:BSDoverQ5}, we will first determine the endomorphism ring of $E$.
\begin{lemma}\label{lemma:EhasCM}
The elliptic curve $E$ geometrically has complex multiplication by $\Z[\sqrt{-5}]$.
\end{lemma}
\begin{proof}
The Hilbert class polynomial for discriminant -20 is $$x^2 - 1264000 \cdot x - 681472000,$$ see for example \cite[Table 2, p.\ 400]{HilbertTable}.
Its zeros are $632000 \pm 282880 \sqrt{5}$. The $j$-invariant for $E$ is $632000 - 282880 \sqrt{5}$, which proves that $E$ geometrically has complex multiplication by $\Z[\sqrt{-5}]$.  
\end{proof}}}

The following propositions will be used to prove Theorem \ref{thm:BSDoverQ5}.

\begin{proposition}\label{thm:twoexplicitmaps}
Let $K = \Qfive$ and
$$\varphi \colon H_{K} \rightarrow E \colon (x : y : 1) \mapsto \left( \varphi_x : \varphi_y : 1 \right), \textrm{ \, with}$$
$$\varphi_x = \frac{\sqrt{5} \cdot x^2 - \sqrt[4]{5} \cdot x + 1}{x}, \qquad \varphi_y = \frac{-\sqrt[4]{5}^3 \cdot xy + \sqrt{5} \cdot y}{x^2}$$
Then the map $\psi \colon H_{\Q(\sqrt5\,)} \rightarrow W := \Res^{K}_{\Q(\sqrt5\,)} E$ naturally induced by $\varphi$ induces an isogeny $\nu \colon \Jac H_{\Q(\sqrt5\,)} \rightarrow W$ over $\Q(\sqrt5\,)$.
\end{proposition}

\begin{proof}
For the Weil restriction we have $$W_{K} = E \times E',$$ where $E'$ over $K$ is the pull-back of $E$ under the automorphism $\sigma\colon \sqrt[4]{5} \mapsto -\sqrt[4]{5}$ of $K$ over $\Q(\sqrt{5}\,)$. Using this identification, after base change, the map $\psi$ becomes $$\psi_K : H_{K} \stackrel{(\varphi, \varphi^\sigma)}\longrightarrow E \times E'.$$

Suppose that the map $\nu_{K}$ induced by $\psi_{K}$ is not an isogeny. Then the image of $\nu_{K}$ in $E \times E'$ is an elliptic curve $F$ over $K$ and we have the following diagram.
$$\xymatrix{
											&	&								&	&						&E		\\ \\
H_K \ar[rr]^{\qquad \geq 2}\ar[rrrrruu]^{\varphi, 2}\ar[rrrrrdd]_{\varphi^\sigma, 2}	&	&F \ar[rr]\ar[rrruu]_1\ar[rrrdd]^1	&	&E \times E' \ar[ruu]\ar[rdd]	&		\\ \\
											&	&								&	&						&E'		}$$

As the morphisms $\varphi$ and $\varphi^\sigma$ are of degree 2, and the morphism $H_K \rightarrow F = \nu(H_K)$ is of degree at least 2, the two morphisms $F \rightarrow E$ and $F \rightarrow E'$ are of degree 1 and defined over $K$. Hence, $E$ and $E'$ must be isomorphic over $K$. Even though $E$ and $E'$ are isomorphic over $\Qfivei$, it is easily verified that they are not isomorphic over $K$. Therefore, $\nu_K$ must be an isogeny and hence also $\nu$ is an isogeny.
\end{proof}

\begin{remark}
The map $\varphi \colon H_K \rightarrow E$ is the quotient of $H_K$ by the automorphism $$H_K \rightarrow H_K \colon \quad x \mapsto \frac{1}{\sqrt5 \cdot x}, \quad y \mapsto \frac{-y}{\sqrt[4]{5}^3 \cdot x^3}.$$
In fact, the geometric automorphism group of $H$ is the dihedral group $D_4$ of order 8, and the Jacobian of any curve of genus 2 over $\Q$ whose automorphism group is non-abelian, is isogenous to the square of an elliptic curve, over a finite extension of $\Q$, cf.\ \cite[Lem.\ 2.4, p.\ 42]{CardonaGonzalezLarioRio}. Remark that this result does not give control on the degree of the field extension needed to define the isogeny.
\end{remark}

Now let us generalise the notion of quadratic twists of elliptic curves to abelian varieties over number fields.

\begin{definition}
Let $A$ be an abelian variety over a number field $K$, and let $K \subset L$ be an extension of degree 2. Then the {\em $L$-quadratic twist} of $A$ over $L$ is the twist of $A$ corresponding to the cocycle $\mathrm{Gal}(L/K) \rightarrow \mathrm{Aut}_L(A)$ mapping the non-trivial element $\sigma \in \mathrm{Gal}(L/K)$ to the automorphism $-1 \colon A \rightarrow A$.
\end{definition}

\onlyinsubfile{The following proposition is probably well-known to the experts.
\begin{proposition}
Let $A$ be an abelian variety over a number field $K$, and let $K \subset L$ be an extension of degree 2. Then the Weil restriction $W := \Res_K^L A_L$ of the base change $A_L$ to $K$ is isogenous to the product $A \times A '$, where $A'$ over $K$ is the $L$-quadratic twist of $A$.
\end{proposition}
\begin{proof}
Recall that $\Hom_{K} (T, W) = \Hom_L (T_L, A_L)$ for any scheme $T$ over $K$.
Consider the morphism $\nu \colon A \times A' \rightarrow W$, given by the morphism $$\kappa \colon A_L \times A'_L \rightarrow A_L \colon (x,y) \mapsto x + \rho(y),$$ where the isomorphism $\rho : A_L' \cong A_L$ comes from the twist data. Then the map $$\nu_L \colon A_L \times A'_L \rightarrow A_L \times A_L^{\sigma}$$ is given by $\kappa$ on the first component and $\kappa^{\sigma}$ on the second component, where $\sigma \colon L \rightarrow L$ is the non-trivial element of $\mathrm{Gal}(L/K)$. Then $\kappa^{\sigma}$ is $$A_L \times A_L' \rightarrow A_L^{\sigma} = A_L \colon  (x,y) \mapsto \sigma(x) + \rho(\sigma(y)).$$ As $\rho(\sigma(y)) = -\sigma(\rho(y))$, by definition of the $L$-quadratic twist, we now find that the kernel of $\nu_L$ is finite and that $\nu$ is an isogeny.
\end{proof}

\begin{example}
For example, for an elliptic curve $E \colon y^2 = f(x)$ over $K$, the $L$-quadratic twist is $E' \colon dy^2 = f(x)$ over $K$ and the isomorphism $\rho$ is given by $E_L \rightarrow E'_L \colon (x,y) \mapsto (x, y /\sqrt{d})$, and $$\nu_L \colon E_L \times E_L' \rightarrow E_L \times E_L \colon (x,y) \mapsto (x+\rho(y), \sigma(x-\rho(y))).$$ The kernel of $\nu_L$ consists of the pairs $(x, \rho^{-1}(x))$ where $x \in E[2]$. Hence, the isogeny $E \times E' \rightarrow W$ has degree 4.
\end{example}}

\begin{proposition}\label{prop:BSDfacts}
Let $A$ and $B$ be abelian varieties over a number field $K$, let $K \subset L$ be a finite extension of number fields and let $C$ be an abelian variety over $L$. Then\\[-1.2cm]
\begin{itemize}\itemsep0pt
\item[(1)] BSD holds for $A \times B$ over $K$ if and only if it holds for $A$ and $B$ over $K$;
\item[(2)] if $A$ and $B$ are isogenous over $K$, then BSD holds for $A$ over $K$ if and only if it holds for $B$ over $K$;
\item[(3)] BSD holds for the Weil restriction $\Res_K^L C$ over $K$ if and only if it holds for $C$ over $L$;
\item[(4)] if $L/K$ is quadratic, BSD holds for the base change $A_L$ over $L$ if and only if it holds for $A$ over $K$ and its $L$-quadratic twist $A'$ over $K$.
\end{itemize}
\end{proposition}

\begin{proof}
For (1) and (2), see \cite[p.\ 422]{TateBourbaki}. For (3), see \cite{MilneInventiones}. In the case $L/K$ is a quadratic extension, $\Res^L_K A_L$ is isogenous over $K$ to $A \times A'$, where $A'/K$ is the $L$-quadratic twist of $A$, cf.\ \cite[Thm., p.\ 53]{Kida}. Now (4) follows from (1), (2) and (3).
\end{proof}

\begin{proof}[Proof (Theorem \ref{thm:BSDoverQ5})]
By Proposition \ref{prop:BSDfacts} part (4), BSD holds for $\Jac H$ and $\Jac H'$ over $\Q$ if and only if it holds for $\Jac H_{\Q(\sqrt{5}\,)}$ over $\Q(\sqrt{5}\,)$. The latter is isogenous over $\Q(\sqrt{5}\,)$ to $\Res^{\Q(\sqrt[4]{5}\,)}_{\Q(\sqrt5\,)} E$ by Proposition \ref{thm:twoexplicitmaps}. Hence, by parts (2) and (3) of Proposition \ref{prop:BSDfacts}, BSD holds for $\Jac H_{\Q(\sqrt5\,)}$ over $\Q(\sqrt{5}\,)$ if and only if it holds for $E$ over $\Qfive$.
\end{proof}

Using the methods in \cite{NumericalVerification}, we can numerically verify that the Birch and Swinnerton-Dyer conjecture holds for $\Jac H$ and $\Jac H'$ in the following sense. We numerically verified that the analytic and algebraic rank agree, and we computed all terms except for $|\Sh|$, with more than 20 digits precision. Then we used the conjectural formula to predict the order of $\Sh$. This predicted order, $|\Sh_{\mathrm{an}}|$, appears to equal 1 in both cases. This gives strong evidence for the conjecture, especially since 1 is the square of an integer, which is to be expected according to \cite{PoonenStoll}.

In fact, we found the following values for the BSD-invariants:\\[-1cm]
\begin{center}
\begin{tabular}{c|cc}
													&$\Jac H$						&$\Jac H'$						\\\hline
$r$													&1							&1							\\
$\lim_{s \rightarrow 1} (s-1)^{-r}L(s) $		&4.54183774632835249986 	&4.54183774632835249986 	\\
$R$													&4.70213971014416647713		&0.94042794202883329543		\\
$\Omega$											&1.93181743899697988452		&9.65908719498489942260		\\
$c_{\mathfrak{p}}$												&$c_2 = 1$, $c_5 = 2$			&$c_2 = 1$, $c_5 = 2$			\\
$|J_{\mathrm{tors}}|$								&2							&2							\\
$\Sh_{\mathrm{an}}$									&1.00000000000000000000		&1.00000000000000000000		\\
\end{tabular}
\end{center}

\begin{remark}
The values of these invariants suggest that $\Jac H$ and $\Jac H'$ are isogenous; they all seem to differ by an integer multiple. Since, the numerical verification succeeded for both curved, the author did not try to actually find an isogeny.
\end{remark}

\onlyinsubfile{}

\section{Methodology}

In this section, I will try to answer the question how you find an elliptic curve $E$ over a number field $K$, with $L \subset K$ of degree 2, such that its Weil restriction to $L$ is isogenous to the Jacobian of a hyperelliptic curve of genus 2 defined over $\Q$ as abelian varieties (without fixed polarisation).

\subsection{Which elliptic curves?}

The product of two elliptic curves over a number field, $E$ and $E'$, taken with the associated product polarisation, does not lie in the Jacobi locus in the moduli space of polarised abelian varieties, cf.\ \cite[Satz 2, p.\ 37]{WeilTorelli}. However, in some cases it might happen that the abelian variety has another polarisation which makes it into the Jacobian of a smooth curve of genus 2. Heuristically, most polarised abelian varieties lie in the Jacobi locus, but also most polarised abelian varieties have only one polarisation, up to multiplication by an integer. So, heuristically it is not so clear whether such $E$ and $E'$ actually exist. Hence, we should be looking for elliptic curves $E$ and $E'$, such that $E \times E'$ contains a smooth curve of genus 2.

The work of Hayashida and Nishi, \cite{HayashidaNishi}, contains sufficient conditions on $E$ and $E'$ for this situation to arise. In particular, \cite[Thm., \S 4, p.\ 14]{HayashidaNishi} states: if $E$ and $E'$ have complex multiplication by the principal order of the imaginary quadratic field $\Q(\sqrt{-m})$ and $m$ is not 1, 3, 7 or 15, then $E \times E'$ contains a smooth curve of genus 2.

\subsection{Reconstruction of the hyperelliptic curve}

Assume that $E$ over $K$ geometrically has complex multiplication by $\O_{-m} = \Z[\alpha_m]$, where $$\alpha_m = \begin{cases}
\sqrt{-m} &\textrm{if } m \not\equiv 3 \mod 4;\\
\tfrac12(\sqrt{-m}+1) &\textrm{if } m \equiv 3 \mod 4. 
\end{cases}$$
Now consider the complexification $E_\C$ and fix an embedding of $\O_{-m}$ in $\C$. Then $E_{\C} \cong \C / \Lambda$, where $\Lambda$ is a lattice of the form $\Z \cdot 1 + \Z \cdot \frac{\beta}{\gamma}$ with $\beta$ and $\gamma \neq 0$ generating, as $\Z$-module, an ideal of $\O_{-m}$. Moreover, $E_{\C}$ has a Hermitian form, whose imaginary part, without loss of generality, gives the standard antisymmetric form $$\begin{pmatrix} 0 &1\\ -1 &0 \end{pmatrix}$$ on $\Gamma$, with respect to the basis just given.

The idea is now to consider the complex lattice $\Z\left(\begin{smallmatrix} 1\\ 0\end{smallmatrix}\right) + \Z\left(\begin{smallmatrix} 0\\ 1\end{smallmatrix}\right) + \Z\left(\begin{smallmatrix} \alpha_m\\ 0\end{smallmatrix}\right) + \Z\left(\begin{smallmatrix} 0\\ \alpha_m\end{smallmatrix}\right)$ inside $\C^2$. We try to put other antisymmetric forms on the lattice, and for each such a form, we choose a basis, such that the antisymmetric form with respect to this basis is of the standard form $$\begin{pmatrix} 0 &0 &1 &0 \\ 0 &0 &0 &1 \\ -1 &0 &0 &0 \\ 0 &-1 &0 &0\end{pmatrix}.$$ After this, we apply a transformation in $\mathrm{GL}_2(\C)$ to obtain a basis that is of the form $\left(\begin{smallmatrix} 1\\ 0\end{smallmatrix}\right), \left(\begin{smallmatrix} 0\\ 1\end{smallmatrix}\right), \left(\begin{smallmatrix} v_1\\ v_2\end{smallmatrix}\right), \left(\begin{smallmatrix} w_1\\ w_2\end{smallmatrix}\right)$, cf.\ \cite[\S 5]{Schlichenmaier}. If the antisymmetric form satisfies the Riemann relations, cf.\ \cite[Lem.\ 1.1 \& 1.2, Chap.\ VII, \S 1, p.\ 132]{LangIntro}, then the matrix $$M = \begin{pmatrix} v_1 &w_1\\ v_2 &w_2\end{pmatrix}$$ will be symmetric and its imaginary part will be positive definite, i.e.\ $M$ has the potential to be the small period matrix of a hyperelliptic curve $H$ of genus 2.

One can then evaluate the theta functions in $M$ and use these to reconstruct the Igusa invariants of $H$. These Igusa invariants can only be computed numerically, up to a certain precision, but we expect them to be rational. If the precision is high enough, we can guess the rational values for the Igusa invariants. Then we can use Mestre's algorithm (\cite{Mestre}) to construct a hyperelliptic curve with these Igusa invariants. This part of the reconstruction procedure is explained in more detail in \cite{Weng}.

\subsection{Constructing algebraic maps}

Now we are in the situation that we found an elliptic curve $E$ over $K$ and a hyperelliptic curve $H$ over $\Q$, such that the base change of $E \times E$ and $J := \mathrm{Jac}(H)$ to $\C$ numerically seem to be isogenous. If such an isogeny exists, we know by GAGA that it is algebraisable and defined over a finite extension of $K$. The only problem that remains is to find such an algebraic isogeny explicitly.

It is possible to numerically construct an analytic isogeny $\tau: H_{\C} \rightarrow J_{\C} \rightarrow E_{\C} \times E_{\C}$. We consider two composite maps $$\tau_1, \tau_2: \xymatrix{H_{\C} \ar[r] &E_{\C} \times E_{\C} \ar@<-.5ex>[r]\ar@<.5ex>[r] &E_{\C}}$$ and try to `guess' them. We assume that the map $\tau_1 : H_{\C} \rightarrow E_{\C}$ is of the shape $$(x,y) \mapsto \frac{\sum_{i=0}^{N} \sum_{j=0}^1 a_{i,j} x^i y^j}{\sum_{i=0}^M \sum_{j=0}^1 b_{i,j} x^i y^j},$$ for certain $a_{i,j}, b_{i,j} \in \C$ and $N,M \in \Z_{\geq 0}$. We pick $R := 2N+2M$ complex-valued points $P_k := (\alpha_k, \beta_k) \in H_{\C}(\C)$ for $k = 1, \ldots, R$ and numerically compute $Q_k := \tau_1(P_k)$. Each such point gives rise to a linear equation $$\sum_{i=0}^N\sum_{j=0}^1 a_{i,j} \alpha_k^i\beta_k^j - Q_k \cdot \sum_{i=0}^M \sum_{j=0}^1 b_{i,j} \alpha_k^i \beta_k^j = 0$$ in the coefficients $a_{i,j}$ and $b_{i,j}$. Or, to phrase it in other words, the vector of coefficients $(a_{0,0}, \ldots, a_{N,1}, b_{0,0}, \ldots, b_{M,1})$ is in the kernel of the matrix $$A = \begin{pmatrix}\alpha_1^0 \beta_1^0 &\cdots &\alpha_1^N \beta_1^1 &-Q_1 \alpha_1^0 \beta_1^0 & \cdots &-Q_1 \alpha_1^M \beta_1^0\\ \vdots &\ddots &\vdots &\vdots &\ddots &\vdots \\ \alpha_{R}^0 \beta_{R}^0 & \cdots &\alpha_{R}^N \beta_{R}^1 &-Q_{R}\alpha_{R}^0 \beta_{R}^0 &\cdots &-Q_{R}\alpha_{R}^M \beta_{R}^1  \end{pmatrix}.$$
We can compute this kernel numerically and choose $N$ and $M$ such that the kernel is 1-dimensional. In this way, we can be sure to find a basis vector, which is a $\C$-multiple of a vector with algebraic entries, instead of obtaining a random $\C$-linear combination of two or more.

We compute a generator for the kernel and rescale it to make one of the non-zero entries equal to 1. Then we use LLL to guess algebraic relations for the other entries. In this way, we found a solution $(a_{0,0}, \ldots, b_{M,1}) \in \overline{\Q}^R$ and, if $M$ and $N$ were chosen appropriately, it can be verified algebraically that these functions indeed define a morphism $\varphi : H_L \rightarrow E_L$, where $L = K(a_{0,0}, \ldots, b_{M,1})$, whose base change to $\C$ is $\tau_1$.

\subsection{Smaller fields}

A priori, the field $L$ might be way too big for a feasible numerical verification of BSD. For example, in our specific case, a priori the curve $H$ and $E$ were defined over $\Q$ and $\Q(\sqrt5\,)$, respectively, but the maps $\varphi$ and $\psi$ were only defined over $L = \Q(\sqrt[8]{5}, i)$ and $\varphi \colon H \rightarrow E \colon (x:y:1) \mapsto (\varphi_x : \varphi_y : 1)$ was given by
$$\varphi_x = \frac{\tfrac12  i \sqrt[4]{5} \cdot x^4 - x^3 - \tfrac12 i \left(\tfrac45 \sqrt[4]{5}^3 - \sqrt[4]{5}\right) \cdot x^2 + \tfrac15\sqrt5 \cdot x + \tfrac1{10}  i  \sqrt[4]{5}}{x^3 + \frac{2i}{5}\sqrt[4]{5}^3 \cdot x^2 - \tfrac15 \sqrt5 \cdot x },$$
$$\varphi_y = \frac{\tfrac14 \varepsilon \sqrt[8]{5}^3 \cdot x^4y + \delta\sqrt[8]{5} \cdot x^3y - \tfrac14 \varepsilon \left( \tfrac45 \sqrt[8]{5}^7 + \sqrt[8]{5}^3 \right) \cdot x^2y - \tfrac{\delta}5\sqrt[8]{5}^{5} \cdot xy + \tfrac1{20} \varepsilon  \sqrt[8]{5}^{3} \cdot y }{x^5 + \frac{3i}{5}\sqrt[4]{5}^3 \cdot x^4 - \tfrac35 \sqrt{5} \cdot x^3 - \tfrac15i\sqrt[4]{5} \cdot x^2},$$
where $\varepsilon = 1-i$ and $\delta = 1+i$. Of course this still proves that $\Jac{H}_L$ and $E_L \times E_L$ are isogenous. 

However, it is not feasible yet to numerically verify BSD for $H_L$. The situation is not as good as in Proposition \ref{prop:BSDfacts} part (4). In the isogeny decomposition of the Weil restriction $\Res^L_{\Q(\sqrt5)}{\Jac(H_L)}$, there will not only be twists of $\Jac H$ occuring, but also higher dimensional factors, see also\ \cite{DiemNaumann}. Even if we are lucky, and all these factors are Jacobians of hyperelliptic curves over $\Q$, these curves will be of genus greater than 3. Numerical verification of BSD for such curves might take too much time.

In order to reduce the size of $L$ and reduce to the case of a quadratic extensions of field, we performed some twists, for example on $E$ by $\varepsilon \sqrt[8]{5}$ and on $H$ by $-1$. We then repeated the procedure in the previous paragraph and even managed to find a map of smaller degree over the smaller field $\Q(\sqrt[4]{5})$.

Having found the appropriate map defined over $\Q(\sqrt[4]{5})$, we were able to get the result in Proposition \ref{thm:twoexplicitmaps} in order to finally prove Theorem \ref{thm:BSDoverQ5}.

\onlyinsubfile{}

\end{document}